# STOCHASTIC PARTIAL DIFFERENTIAL EQUATIONS DRIVEN BY LÉVY SPACE-TIME WHITE NOISE

By Arne Løkka,[1] Bernt Øksendal and Frank Proske

*University of Oslo, Norwegian School of Economics and Business Administration and University of Oslo, and University of Oslo*

In this paper we develop a white noise framework for the study of stochastic partial differential equations driven by a $d$-parameter (pure jump) Lévy white noise. As an example we use this theory to solve the stochastic Poisson equation with respect to Lévy white noise for any dimension $d$. The solution is a stochastic distribution process given explicitly. We also show that if $d \leq 3$, then this solution can be represented as a classical random field in $L^2(\mu)$, where $\mu$ is the probability law of the Lévy process. The starting point of our theory is a chaos expansion in terms of generalized Charlier polynomials. Based on this expansion we define Kondratiev spaces and the Lévy Hermite transform.

**1. Introduction.** White noise analysis has become a subject of much current interest. This theory was first treated by Hida [14] and extensively studied in many other works. See [16] and the references therein. These investigations are based on the concept of a Gaussian measure and the associated expansion into Hermite polynomials. Later on an extension of white noise theory to non-Gaussian analysis was established in [4] and developed further by Kondratiev, Da Silva, Streit and Us [24] and Kondratiev, Da Silva and Streit [23]. The main tool of this theory is a biorthogonal decomposition, which extends the Wiener–Itô chaos expansion. White noise analysis has been used in a broad range of applications. This approach was originally applied in quantum physics. See, for example, [3] or [2]. Subsequently, new applications have been found in stochastic (partial) differential equations

Received November 2002; revised February 2003.
[1]Supported by the Norwegian Research Council, Grant 134228/432.

[17]. See also [21] and [6] to mention a few. More recently, the theory has been applied to finance [1]. See [18] and [12] for the fractional Brownian motion case and [10] and [34] in the non-Gaussian case.

The object of this paper is to provide a white noise framework, based on results in [28, 10, 34] and [17], to study SPDEs driven by (pure jump) Lévy processes. We apply this theory to solve the *stochastic Poisson equation driven by a d-parameter ( pure jump) Lévy white noise.* That is, consider the following model for the temperature $U(x)$ at point $x$ in a bounded domain $D$ in $\mathbb{R}^d$. Suppose that the temperature at the boundary $\partial D$ of $D$ is kept equal to zero and that there is a random heat source in $D$ modeled by *Lévy white noise* $\dot{\eta}(x) = \dot{\eta}(x_1, \ldots, x_d)$. Then $U$ is described by the equation

$$
(1.1) \qquad \begin{aligned} \Delta U(x) &= -\dot{\eta}(x), & x = (x_1, \ldots, x_d) \in D \\ U(x) &= 0, & x \in \partial D. \end{aligned}
$$

It is natural to guess that the solution must be

$$
(1.2) \qquad U(x) = U(x, \omega) = \int_D G(x, y) \, d\eta(y),
$$

where $G(x, y)$ is the classical Green function for $D$ and the integral on the right-hand side is a multiparameter Itô integral with respect to the $d$-parameter Lévy process $\eta(x)$. But the integral on the right-hand side of (1.2) only makes sense if $G(x, \cdot)$ is square integrable in $D$ with respect to the Lebesgue measure. The latter is true if and only if the dimension $d$ is chosen lower than 4. Despite this difficulty we will show the existence of a unique explicit solution

$$
x \mapsto U(x, \cdot) \in (\mathcal{S})_{-1},
$$

where $(\mathcal{S})_{-1}$ is a suitable space of stochastic distributions, called the Kondratiev space.

The stochastic Poisson equation (1.1) was discussed by Walsh [39] in the case of Brownian white noise $W$. He proved that there exists for all $d$ a Sobolev space $H^{-n}(D)$ and an $H^{-n}(D)$-valued stochastic process

$$
U = U(\omega) : \Omega \to H^{-n}(D)
$$

such that (1.1) holds in the sense of distributions, for example,

$$
\langle U(\cdot, \omega), \Delta \phi \rangle = -\langle W(\cdot, \omega), \phi \rangle \qquad \text{a.s. for all } \phi \in H^{-n}(D).
$$

The solution of Walsh is given explicitly by

$$
(1.3) \qquad \langle U, \phi \rangle = \int_{\mathbb{R}^d} \int_{\mathbb{R}^d} G(x, y) \phi(x) \, dx \, dB(y); \qquad \phi \in H^{-n}(D).
$$

The system (1.1) was also studied in [17] in the Gaussian case. There the solution $U(x)$, which takes values in the Kondratiev space, can be described by its action on test functions $f \in (\mathcal{S})_1$,

$$(1.4) \qquad \langle U(x), f \rangle = \int_{\mathbb{R}^d} G(x,y) \langle W(y), f \rangle \, dy; \qquad f \in (\mathcal{S})_1.$$

If we compare (1.3) and (1.4) we find that the Walsh solution takes $x$-averages for almost all $\omega$, whereas the last one takes $\omega$-averages for all $x$.

Our solution is an extension of (1.4) to Lévy processes. The approach we use to solve (1.1) is based on a chaos expansion in terms of generalized Charlier polynomials (cf. [28]) and on concepts developed in [17, 10] and [34]. Our method, which can be applied to other classes of SPDEs, has the advantage that SPDEs can be interpreted in the usual strong sense with respect to time and space. There is no need for a weak distribution interpretation with respect to time and space. Furthermore, the Walsh construction reveals the disadvantage of defining a multiplication of (Sobolev or Schwartz) distributions, if one considers SPDEs, where the noise is involved multiplicatively. However, on the Kondratiev space $(\mathcal{S})_{-1}$ we can define a multiplication, the *Lévy Wick product*. This gives a natural interpretation of SPDEs, where the noise or other terms appear multiplicatively. Furthermore, in some cases solutions can be explicitly obtained in terms of the Wick product. See [17].

The general machinery, developed in this paper, is of independent interest and we are convinced that it serves a useful tool for the study of a large class of stochastic partial differential equations driven by Lévy space-time white noise.

Finally, let us mention that there has recently been an increasing interest in solving SPDEs driven by $d$-parameter Lévy processes. We refer to [5, 29] and the references therein.

We shall give an overview of the paper. In Section 2 we introduce a white noise framework for the study of SPDEs driven by $d$-parameter Lévy processes. The starting point of our theory is a chaos expansion in terms of generalized Charlier polynomials. Based on this expansion we define Kondratiev spaces, the Wick product and the $d$-parameter Lévy white noise. Further, we give the definition of the Lévy Hermite transform and state a characterization theorem for the Kondratiev space $(\mathcal{S})_{-1}$. In Section 3 we use the tools developed in Section 2 to apply it to solve the stochastic Poisson equation driven by a $d$-parameter Lévy white noise. Finally, in Section 4 we discuss the solution and its properties. In particular, we show that if $d \leq 3$, the solution can be represented as a classical random field in $L^2(\mu)$, where $\mu$ is the underlying probability law of the Lévy process.

**2. Framework.** In this section we give the general framework to be used later. The starting point for our discussion are white noise concepts for Lévy processes, developed in [10, 34] and [28]. Actually, we empasize the use of multidimensional structures, that is, the white noise we intend to consider is indexed by a multidimensional parameter set. Our presentation and notation will follow that of [17] closely, where Gaussian white noise theory is treated. For more information about white noise theory we refer to [16, 26] and [32].

2.1. *A white noise construction of Poisson random measures associated with a Lévy process.* In this paper we confine ourselves to ($d$-parameter) pure jump Lévy processes without drift.

A pure jump Lévy process $\eta(t)$ on $\mathbb{R}$ with no drift is a process with independent and stationary increments, continuous in probability and with no Brownian motion part. The characteristic function of such a process is given by the Lévy–Khintchine formula in terms of the Lévy measure $\nu$ of the Lévy process, that is, in terms of a measure $\nu$ on $\mathbb{R}_0 := \mathbb{R} - \{0\}$, that integrates the function $1 \wedge z^2$. Hence, driftless pure jump Lévy processes can be characterized as Lévy processes with characteristic triplet $(0, 0, \nu)$. For general information about Lévy processes see [8] and [36]. In general, such processes do not possess the chaotic representation property, but they admit a chaos representation with respect to Poisson random measures (see [19]). Therefore, we aim at viewing these processes as elements of a certain Poisson space. In this framework we will give a white noise construction of Poisson random measures and, since our emphasis lies on processes indexed by multidimensional sets, we will define $d$-parameter (pure jump) Lévy processes. Further, we prove a chaos expansion in terms of generalized Charlier polynomials.

A usual starting point in white noise analysis is the application of the Bochner–Minlos theorem to prove the existence of a probability measure on the space of tempered distributions $\mathcal{S}'(\mathbb{R}^d)$. However, it turns out that $\mathcal{S}'(\mathbb{R}^d)$ is not the most appropriate for dealing with Lévy processes since this choice would require restrictive conditions to be imposed on the Lévy measure. This circumstance comes from the fact that the Lévy measure has a singularity at zero. Therefore, we use the construction of a nuclear algebra $\widetilde{\mathcal{S}}(X)$, which is more tractable for our purpose. In fact, the space $\widetilde{\mathcal{S}}(X)$ is a variant of the Schwartz space on $X = \mathbb{R}^d \times \mathbb{R}_0$, more precisely $\widetilde{\mathcal{S}}(X)$ is a subspace of the Schwartz space modulo a certain subspace depending on the Lévy measure. Let us first give the construction of $\widetilde{\mathcal{S}}(X)$ (cf. [28]).

In the following let $\{\xi_n\}_{n \geq 0}$ be the complete orthogonal system of $L^2(\mathbb{R})$, consisting of the *Hermite functions.* Then the (countably Hilbertian) nuclear topology of the Schwartz space $\mathcal{S}(\mathbb{R}^d)$ is induced by the compatible system of norms

$$\|\varphi\|_\gamma^2 := \sum_{\alpha \in \mathbb{N}^d} (1+\alpha)^{2\gamma} (\varphi, \xi_\alpha)^2_{L^2(\mathbb{R}^d)}, \qquad \gamma \in \mathbb{N}_0^d, \tag{2.1}$$

where $\xi_\alpha := \prod_{i=1}^d \xi_{\alpha_i}$ and $(1+\alpha)^{2\gamma} := \prod_{i=1}^d (1+\alpha_i)^{2\gamma_i}$ for $\alpha = (\alpha_1, \ldots, \alpha_d) \in \mathbb{N}^d$ and $\gamma = (\gamma_1, \ldots, \gamma_d) \in \mathbb{N}_0^d$. Now let us take a numbering of the norms in (2.1), say $\|\cdot\|_{\gamma_i}$, and define the norms $\|\varphi\|_p = \sum_{i=1}^p \|\varphi\|_{\gamma_i}$, $p \in \mathbb{N}$. Then $\|\cdot\|_p$, $p \in \mathbb{N}$ are increasing pre-Hilbertian norms on $\mathcal{S}(\mathbb{R}^d)$. It is well known that the norms $\|\cdot\|_p$ are equivalent to the norms $\|\cdot\|_{q,\infty}$, given by

$$(2.2) \quad \|\varphi\|_{q,\infty} := \sup_{0 \le k, |\gamma| \le q} \sup_{z \in \mathbb{R}^d} |(1+|z|^k)\partial^\gamma \varphi(z)|, \quad q \in \mathbb{N}_0,$$

where $\partial^\gamma \varphi = \frac{\partial^{|\gamma|}}{\partial z_1^{\gamma_1} \cdots \partial z_d^{\gamma_d}} \varphi$ for $\gamma = (\gamma_1, \ldots, \gamma_d) \in \mathbb{N}_0^d$ with $|\gamma| := \gamma_1 + \cdots + \gamma_d$. We mention the following important property of the norms $\|\cdot\|_p$ (see [20]): For all $p \in \mathbb{N}$ there exists a constant $M_p > 0$ such that for all $\varphi, \psi \in \mathcal{S}(\mathbb{R}^d)$,

$$(2.3) \quad \|\varphi\psi\|_p \le M_p \|\varphi\|_p \|\psi\|_p.$$

We then define the space $\mathcal{S}(X)$ by

$$(2.4)\ \mathcal{S}(X) := \left\{ \varphi \in \mathcal{S}(\mathbb{R}^{d+1}) : \varphi(z_1, \ldots, z_d, 0) = \left(\frac{\partial}{\partial z_{d+1}} \varphi\right)(z_1, \ldots, z_d, 0) = 0 \right\}.$$

It follows that $\mathcal{S}(X)$ is a closed subspace of $\mathcal{S}(\mathbb{R}^{d+1})$. Thus, $\mathcal{S}(X)$ is a (countably Hilbertian) nuclear space with respect to the restriction of the norms $\|\cdot\|_p$. Moreover, it is a nuclear algebra, that is, $\mathcal{S}(X)$ is, in addition, a topological algebra with respect to the multiplication of functions. In the sequel we denote by $\lambda^{\times d}$ the Lebesgue measure on $\mathbb{R}^d$ and by $\nu$ a Lévy measure of on $\mathbb{R}_0$. We set $\pi = \lambda^{\times d} \times \nu$. We need the following result.

LEMMA 2.1. *There exists an element denoted by $1 \otimes \dot{\nu}$ in $\mathcal{S}'(X)$ such that*

$$\langle 1 \otimes \dot{\nu}, \phi \rangle = \int_X \phi(y) \pi(dy)$$

*for all $\phi \in \mathcal{S}(X)$, where $\langle 1 \otimes \dot{\nu}, \phi \rangle = (1 \otimes \dot{\nu})(\phi)$ is the action of $1 \otimes \dot{\nu}$ on $\phi$. The notation $\dot{\nu}$ shall indicate that $\dot{\nu}$ is the Radon–Nikodym derivative of $\nu$ in a generalized sense.*

PROOF. Without loss of generality we consider the case $d = 1$. Set $L(\varphi) = \int_X \varphi(z) \pi(dz)$. Let $\varphi_n, \varphi \in \mathcal{S}(X)$ with $\varphi_n \to \varphi$ in $\mathcal{S}(X)$. By Taylor's formula we have for $\varphi \in \mathcal{S}(X)$ that

$$\varphi(x,z) = \varphi(x,0) + \left(\frac{\partial}{\partial z}\varphi\right)(x,0)z + \frac{1}{2}\left(\frac{\partial^2}{\partial^2 z}\varphi\right)(x,\xi)z^2$$

$$= \frac{1}{2}\left(\frac{\partial^2}{\partial^2 z}\varphi\right)(x,\xi)z^2$$

for a point $\xi$ between $0$ and $z$. We assume w.l.o.g. that the measure $\nu$ vanishes outside of $[-1,0) \cup (0,1]$. Therefore, it follows by (2.2) that

$$\begin{aligned}|L(\varphi_n - \varphi)| &\leq \int_{\mathbb{R}} \int_{-1}^{1} |\varphi_n(x,z) - \varphi(x,z)| \nu(dz) \lambda(dx) \\ &\leq \int_{\mathbb{R}} \int_{-1}^{1} \frac{(1+|x|^2+|z|^2)|\varphi_n(x,z) - \varphi(x,z)|}{z^2} \\ &\qquad \times \frac{z^2}{(1+|x|^2)} \nu(dz) \lambda(dx) \\ &\leq \|\varphi_n - \varphi\|_{2,\infty} \int_{\mathbb{R}} \frac{1}{(1+|x|^2)} \lambda(dx) \int_{-1}^{1} z^2 \nu(dz) \\ &\to 0 \qquad \text{for } n \to \infty.\end{aligned}$$

So the linear functional $L$ is continuous on $\mathcal{S}(X)$. $\square$

Next define the space

(2.5) $$\mathcal{N}_\pi := \{\phi \in \mathcal{S}(X) : \|\phi\|_{L^2(\pi)} = 0\}.$$

By the same arguments as in the proof of Lemma 2.1 it can be shown that $\mathcal{N}_\pi$ is a closed subspace of $\mathcal{S}(X)$. Furthermore, one checks that it is a closed ideal of $\mathcal{S}(X)$. Now we introduce the space $\widetilde{\mathcal{S}}(X)$, which we use to construct the white noise measure.

DEFINITION 2.2. We define the space $\widetilde{\mathcal{S}}(X)$ as follows,

(2.6) $$\widetilde{\mathcal{S}}(X) = \mathcal{S}(X)/\mathcal{N}_\pi.$$

The space $\widetilde{\mathcal{S}}(X)$ is a (countably Hilbertian) nuclear algebra with the compatible system of norms

(2.7) $$\|\widehat{\phi}\|_{p,\pi} := \inf_{\psi \in \mathcal{N}_\pi} \|\phi + \psi\|_p, \qquad p \in \mathbb{N}.$$

See page 72 in [13]. Further, let $\widetilde{\mathcal{S}}'(X)$ denote the topological dual of $\widetilde{\mathcal{S}}(X)$.

We obtain the following corollary to Lemma 2.1:

COROLLARY 2.3. *The functional* $L(\widehat{\phi}) := \int_X \phi(z) \pi(dz)$ *satisfies the inequality*

$$|L(\widehat{\phi})| \leq M_p \|\widehat{\phi}\|_{p,\pi}$$

*for all* $p \geq p_0$, *which yields the continuity of the functional* $L$ *on* $\widetilde{\mathcal{S}}(X)$.

THEOREM 2.4. *There exists a unique probability measure $\mu$ on the Borel sets of $\widetilde{\mathcal{S}}'(X)$ with the following Poissonian characteristic functional with intensity $\pi$ such that for all $\phi \in \widetilde{\mathcal{S}}(X)$:*

$$(2.8) \qquad \int_{\widetilde{\mathcal{S}}'(X)} e^{i\langle\omega,\phi\rangle} \, d\mu(\omega) = \exp\left(\int_X (e^{i\phi} - 1) \, d\pi\right),$$

*where $\langle\omega,\phi\rangle = \omega(\phi)$ is the action of $\omega \in \widetilde{\mathcal{S}}'(X)$ on $\phi \in \widetilde{\mathcal{S}}(X)$. Moreover, there exists a $p_0 \in \mathbb{N}$ such that $1 \otimes \dot{\nu} \in \widetilde{\mathcal{S}}_{-p_0}(X)$ and a natural number $q_0 > p_0$ such that the imbedding operator $\widetilde{\mathcal{S}}_{q_0}(X) \hookrightarrow \widetilde{\mathcal{S}}_{p_0}(X)$ is Hilbert–Schmidt and $\mu(\widetilde{\mathcal{S}}_{-q_0}(X)) = 1$. The space $\widetilde{\mathcal{S}}_p(X)$ denotes the completion of $\widetilde{\mathcal{S}}(X)$ with respect to $\|\cdot\|_{p,\pi}$ and $\widetilde{\mathcal{S}}_{-p}(X)$ is the corresponding dual with norm $\|\cdot\|_{-p,\pi}$.*

PROOF. Since $|e^{iz} - 1| \leq |z|$, the result follows from Corollary 2.3 and Bochners theorem for conuclear spaces [13]. □

We call the probability measure $\mu$ on $\Omega = \widetilde{\mathcal{S}}'(X)$ in Theorem 2.4 *Lévy white noise probability measure*. It turns out that this measure satisfies the *first condition of analyticity* in the following sense (see [23]).

LEMMA 2.5. *The Lévy white noise measure $\mu$ satisfies the first condition of analyticity, that is, there exists $\epsilon > 0$ and a $p_0$ such that*

$$\int_{\widetilde{\mathcal{S}}'(X)} \exp(\epsilon \|\omega\|_{-p_0}) \, d\mu(\omega) < \infty.$$

PROOF. The proof follows the argument of Lemma 3 in [38]. Introduce the moment functions of $\mu$, which by a criterion of Cramér [9] can be expressed by

$$M_n(\phi) := \int_{\widetilde{\mathcal{S}}'(X)} \langle\omega,\phi\rangle^n \, d\mu(\omega) = \frac{d}{dt^n} L(t\phi)\bigg|_{t=0}$$

for every $\phi \in \widetilde{\mathcal{S}}(X)$, $n \in \mathbb{N}$. Define the set

$$\Lambda_n^k := \left\{(\alpha_1, \ldots, \alpha_k) \in \mathbb{N}^k : \sum_{i=1}^k \alpha_i = n\right\}.$$

Then we obtain the following expression for $M_n$:

$$(2.9) \qquad M_n(\phi) = \sum_{k=1}^n \frac{n!}{k!} \sum_{\alpha \in \Lambda_n^k} \prod_{j=1}^k \frac{\langle 1 \otimes \dot{\nu}, \phi^{\alpha_j}\rangle}{\alpha_j!}.$$

We get for the number $p_0$ in Theorem 2.4 that

$$|\langle 1 \otimes \dot{\nu}, \phi\rangle| \leq \|1 \otimes \dot{\nu}\|_{-p_0,\pi} \|\phi\|_{p_0,\pi} < \infty.$$

Next relation (2.3) implies that for all $p \in \mathbb{N}$ there exists a constant $M_p > 0$ such that for all $\phi, \psi \in \widetilde{\mathcal{S}}(X)$,

(2.10) $$\|\phi\psi\|_{p,\pi} \leq M_p \|\phi\|_{p,\pi} \|\psi\|_{p,\pi}.$$

Thus, we get that

$$|\langle 1 \otimes \dot{\nu}, \phi^{\alpha_j}\rangle| \leq \|1 \otimes \dot{\nu}\|_{-p_0,\pi} (M_{p_0})^{\alpha_j} \|\phi\|_{p_0,\pi}^{\alpha_j},$$

if we choose $M_{p_0} \geq 1$. So we deduce from (2.9) that

$$|M_n(\phi)| \leq \sum_{k=1}^{n} \frac{n!}{k!} \sum_{\alpha \in \Lambda_n^k} \prod_{j=1}^{k} \frac{\|1 \otimes \dot{\nu}\|_{-p_0,\pi}}{\alpha_j!} C_{p_0}^n \|\phi\|_{p_0,\pi}^n$$
$$= F_n(\|1 \otimes \dot{\nu}\|_{-p_0,\pi}) C_{p_0}^n \|\phi\|_{p_0,\pi}^n,$$

where $F_n(x)$ is the $n$th moment of the Poisson distribution with intensity $x$ and where $C_{p_0}$ is a constant. Further, it is known that for a Poisson distribution with intensity $x = \|1 \otimes \dot{\nu}\|_{-p_0,\pi}$, there exists a constant $C_x$ such that for all $n \in \mathbb{N}$,

$$|F_n(\|1 \otimes \dot{\nu}\|_{-p_0,\pi})| \leq n! C_{\|1 \otimes \dot{\nu}\|_{-p_0,\pi}}^n.$$

Therefore, we get for a $C > 0$ that

$$|M_n(\phi)| \leq n! C^n \|\phi\|_{p_0,\pi}^n.$$

The claimed result follows from Lemma 3 in [23]. □

Further, consider the function $\alpha$ defined by $\alpha(\phi) = \log(1 + \varphi) \mod \mathcal{N}_\pi$ for $\phi = \widehat{\varphi}$ with $\varphi(x) > -1$. Note that $\alpha$ is holomorphic at zero and invertible. With the help of Lemma 2.5, it can be shown just as in [28] that there exist symmetric kernels $C_n(\omega)$ such that for all $\phi$ in an open neighborhood of zero in $\widetilde{\mathcal{S}}(X)$,

(2.11) $$\widetilde{e}(\phi, \omega) := \frac{\exp\langle\omega, \alpha(\phi)\rangle}{E_\mu[e^{\langle\omega,\alpha(\phi)\rangle}]} = \sum_{n\geq 0} \frac{1}{n!} \langle C_n(\omega), \phi^{\otimes n}\rangle,$$

where $\phi^{\otimes n} \in \widetilde{\mathcal{S}}(X)^{\widehat{\otimes} n}$. The symbol $\widetilde{\mathcal{S}}(X)^{\widehat{\otimes} n}$ denotes the $n$th completed symmetric tensor product of $\widetilde{\mathcal{S}}(X)$ with itself. The elements of this space can be seen as functions $f \in \mathcal{S}(X^n)$ modulo $\mathcal{N}_{\pi \times n}$ such that $f = f(x_1, \ldots, x_n)$ is symmetric with respect to the variables $x_1, \ldots, x_n \in X$. From (2.11) we conclude that the $C_n$ are *generalized Charlier polynomials* (see [23]). We have that

(2.12) $$\{\langle C_n(\omega), \phi^{(n)}\rangle : \phi^{(n)} \in \widetilde{\mathcal{S}}(X)^{\widehat{\otimes} n}, n \in \mathbb{N}_0\}$$

is a total set in $L^2(\mu)$. Furthermore, for all $n$, $m$, $\phi^{(n)} \in \widetilde{\mathcal{S}}(X)^{\widehat{\otimes}n}$ and $\psi^{(m)} \in \widetilde{\mathcal{S}}(X)^{\widehat{\otimes}m}$ the orthogonality relation

$$(2.13) \quad \int_{\widetilde{\mathcal{S}}'(X)} \langle C_n(\omega), \phi^{(n)} \rangle \langle C_m(\omega), \psi^{(m)} \rangle \, d\mu(\omega) = \delta_{n,m} n! (\phi^{(n)}, \psi^{(n)})_{L^2(X^n)}$$

holds. See [28].

REMARK 2.6. It can be easily seen from (2.13) and the construction of $\widetilde{\mathcal{S}}(X)$ that the Lévy white noise measure $\mu$ is *nondegenerate* in the following sense (see [23]): Let $F$ be a continuous polynomial, that is, $F$ is of the form $F(\omega) = \sum_{n=1}^{N} \langle \omega^{\otimes n}, \phi^{(n)} \rangle$ for $\omega \in \widetilde{\mathcal{S}}'(X)$, $N \in \mathbb{N}_0$ with $\phi^{(n)} \in \widetilde{\mathcal{S}}_{\mathbb{C}}(X)^{\widehat{\otimes}n}$ [complexification of $\widetilde{\mathcal{S}}(X)^{\widehat{\otimes}n}$]. If $F = 0$ $\mu$-a.e., then $F(\omega) = 0$ for *all* $\omega \in \widetilde{\mathcal{S}}'(X)$. We mention that this property is essential for the construction of certain test function and distribution spaces (see [23, 28]).

Next, for functions $f : X^n \to \mathbb{R}$ define the *symmetrization* $(f)^\wedge$ of $f$ by

$$(2.14) \quad (f)^\wedge(x_1, \ldots, x_n) := \frac{1}{n!} \sum_n f(x_{\sigma_1}, \ldots, x_{\sigma_n})$$

for all permutations $\sigma$ of $\{1, \ldots, n\}$. Then a function $f : X^n \to \mathbb{R}$ is symmetric, if and only if $\widehat{f} = f$. Denote by $\widehat{L}^2(X^n, \pi^{\times n})$ the space of all symmetric functions on $X^n$, which are square integrable with respect to $\pi^{\times n}$. Let $f_n \in \widehat{L}^2(X^n, \pi^{\times n})$. Since $\mathcal{S}(X)$ is dense in $L^2(X, \pi)$ (cf. [28]), we can choose a sequence $f_n^{(i)}$ in $\widetilde{\mathcal{S}}(X)^{\widehat{\otimes}n}$ with $f_n^{(i)} \to f_n$ in $L^2(X^n, \pi^{\times n})$. Then (2.13) implies the existence of a well defined $\langle C_n(\omega), f_n \rangle$ such that

$$(2.15) \quad \langle C_n(\omega), f_n \rangle = \lim_i \langle C_n(\omega), f_n^{(i)} \rangle \quad \text{in } L^2(X^n, \pi^{\times n}).$$

Since $C_1(\omega) = \omega - 1 \otimes \dot{\nu}$ for all $\phi \in \widetilde{\mathcal{S}}(X)$ (see [28]), we get

$$(2.16) \quad \int_{\widetilde{\mathcal{S}}'(X)} \langle \omega - 1 \otimes \dot{\nu}, f \rangle^2 \, d\mu(\omega) = \|f\|_{L^2(\pi)}^2.$$

Further, if we define for Borelian $\Lambda_1 \subset \mathbb{R}^d$, $\Lambda_2 \subset \mathbb{R}_0$ with $\pi(\Lambda_1 \times \Lambda_2) < \infty$ the random measures

$$(2.17) \quad \begin{aligned} N(\Lambda_1, \Lambda_2) &:= \langle \omega, \chi_{\Lambda_1 \times \Lambda_2} \rangle \quad \text{and} \\ \widetilde{N}(\Lambda_1, \Lambda_2) &:= \langle \omega - 1 \otimes \dot{\nu}, \chi_{\Lambda_1 \times \Lambda_2} \rangle, \end{aligned}$$

we see from their characteristic functions that $N$ is a Poisson random measure and $\widetilde{N}$ is the corresponding compensated Poisson random measure. The compensator of $N(\Lambda_1, \Lambda_2)$ is given by $\pi$. Therefore, it is natural to define the stochastic integral of $\phi \in L^2(\pi)$ with respect to $\widetilde{N}$ by

$$\int_X \phi(x,z) \widetilde{N}(dx, dz) := \langle \omega - 1 \otimes \dot{\nu}, \phi \rangle. \tag{2.18}$$

In particular, if we define

$$\widetilde{\eta}(x) := \int_X \chi_{[0,x_1] \times \cdots \times [0,x_d]}(x) \cdot z \widetilde{N}(dx, dz) \tag{2.19}$$
$$\text{for } x = (x_1, \ldots, x_d) \in \mathbb{R}^d,$$

where $[0, x_i]$ is interpreted as $[x_i, 0]$, if $x_i < 0$ and where the Lévy measure $\nu$ is assumed to integrate $z^2$, then $\widetilde{\eta}(x)$ has a version $\eta(x)$, which is cadlag in each component $x_i$. This follows with the help of (2.13). We call $\eta(x)$ *d-parameter Lévy process* or *space-time Lévy process*.

We conclude this section with a chaos expansion result in terms of the generalized Charlier polynomials $C_n$. The result is a consequence of (2.12) and (2.13).

THEOREM 2.7. *If $F \in L^2(\mu)$, then there exists a unique sequence $f_n \in \widehat{L}^2(X^n)$ such that*

$$F(\omega) = \sum_{n \geq 0} \langle C_n(\omega), f_n \rangle. \tag{2.20}$$

*Moreover, we have the isometry*

$$\|F\|^2_{L^2(\mu)} = \sum_{n \geq 0} n! \|f_n\|^2_{L^2(X^n)}. \tag{2.21}$$

2.2. *Chaos expansion, Kondratiev spaces $(\mathcal{S})_\rho, (\mathcal{S})_{-\rho}$ and Lévy white noise.* First we reformulate the chaos expansion of Theorem 2.7. Then we use the new expansion to define a Wick product on spaces of stochastic test functions and stochastic distributions. The definitions and results here are analogous to the one-parameter case, which is treated in [10] and [34].

From now on we suppose that our Lévy measure $\nu$ satisfies the condition of [30], namely, that for every $\varepsilon > 0$ there exists a $\lambda > 0$ such that

$$\int_{\mathbb{R} \setminus (-\varepsilon, \varepsilon)} \exp(\lambda |z|) \nu(dz) < \infty. \tag{2.22}$$

This implies that our Lévy measure has finite moments of all orders $\geq 2$.

For later use we introduce multi-indices of arbitrary length. To simplify the notation, we regard multi-indices as elements of the space $(\mathbb{N}_0^\mathbb{N})_c$ of all sequences $\alpha = (\alpha_1, \alpha_2, \ldots)$ with elements $\alpha_i \in \mathbb{N}_0$ and with compact support, that is, with only finitely many $\alpha_i \neq 0$. We define

$$\mathcal{J} = (\mathbb{N}_0^\mathbb{N})_c.$$

Further, we set $\mathrm{Index}(\alpha) = \max\{i : \alpha_i \neq 0\}$ and $|\alpha| = \sum_i \alpha_i$ for $\alpha \in \mathcal{J}$.

Next we consider two families of orthogonal polynomials. We use these polynomials to reformulate the chaos expansion of Theorem 2.7. First let $\{\xi_k\}_{k \geq 1}$ be the Hermite functions, just as in Section 2.1. Now choose a bijective map

$$h : \mathbb{N}^d \to \mathbb{N}.$$

Define the function $\zeta_k(x_1, \ldots, x_d) = \xi_{i_1}(x_1) \cdots \xi_{i_d}(x_d)$, if $k = h(i_1, \ldots, i_d)$ for $i_j \in \mathbb{N}$. Then $\{\zeta_k\}_{k \geq 1}$ constitutes an orthonormal basis of $L^2(\mathbb{R}^d)$.

Further, let $\{l_m\}_{m \geq 0}$ be the orthogonalization of $\{1, z, z^2, \ldots\}$ with respect to the innerproduct of $L^2(\varrho)$, where $\varrho(dz) = z^2 \nu(dz)$. Then define the polynomials

$$(2.23) \qquad p_m(z) = \frac{1}{\|l_{m-1}\|_{L^2(\rho)}} z \cdot l_{m-1}(z).$$

The polynomials $p_m$ form a *complete* orthonormal system in $L^2(\nu)$ (see [34]). We shall mention that we could also use any orthonormal basis in $\mathcal{S}(X) \subset L^2(\nu)$ for $d = 0$ instead of the polynomials $p_m$. In this case the integrability condition (2.22) reduces to the requirement of the existence of the second moment with respect to $\nu$. The choice of the polynomials $p_m$ serves to ease notation.

Next define the bijective map

$$(2.24) \qquad z : \mathbb{N} \times \mathbb{N} \to \mathbb{N}; \qquad (i,j) \mapsto j + (i+j-2)(i+j-1)/2.$$

Note that $z(i,j)$ gives the "Cantor diagonalization" of $\mathbb{N} \times \mathbb{N}$.

Then, if $k = z(i,j)$ for $i, j \in \mathbb{N}$, let

$$\delta_k(x,z) = \zeta_i(x) p_j(z).$$

Further, assume $\mathrm{Index}(\alpha) = j$ and $|\alpha| = m$ for $\alpha \in \mathcal{J}$ and identify the function $\delta^{\otimes \alpha}$ as

$$(2.25) \quad \begin{aligned} & \delta^{\otimes \alpha}((x_1, z_1), \ldots, (x_m, z_m)) \\ &= \delta_1^{\otimes \alpha_1} \otimes \cdots \otimes \delta_j^{\otimes \alpha_j}((x_1, z_1), \ldots, (x_m, z_m)) \\ &= \delta_1(x_1, z_1) \cdots \delta_1(x_{\alpha_1}, z_{\alpha_1}) \\ & \quad \cdots \delta_j(x_{\alpha_1 + \cdots + \alpha_{j-1}+1}, z_{\alpha_1 + \cdots + \alpha_{j-1}+1}) \cdots \delta_j(x_m, z_m), \end{aligned}$$

where the terms with zero-components $\alpha_i$ are set equal to 1 in the product ($\delta_i^{\otimes 0} = 1$).

Finally, we define the *symmetrized tensor product* of the $\delta_k$'s, denoted by $\delta^{\widehat{\otimes} \alpha}$ as

$$(2.26) \quad \begin{aligned} & \delta^{\widehat{\otimes} \alpha}((x_1, z_1), \ldots, (x_m, z_m)) \\ &= (\delta^{\otimes \alpha})^{\wedge}((x_1, z_1), \ldots, (x_m, z_m)) \\ &= \delta_1^{\widehat{\otimes} \alpha_1} \widehat{\otimes} \cdots \widehat{\otimes} \delta_j^{\widehat{\otimes} \alpha_j}((x_1, z_1), \ldots, (x_m, z_m)). \end{aligned}$$

For $\alpha \in \mathcal{J}$ define

$$K_\alpha(\omega) := \langle C_{|\alpha|}(\omega), \delta^{\widehat{\otimes}\alpha} \rangle, \tag{2.27}$$

where we let $K_0(\omega) = 1$. For example, if $\alpha = \epsilon^l$ with

$$\epsilon^l(j) = \begin{cases} 1, & \text{for } j = l, \\ 0, & \text{else,} \end{cases} \quad l \geq 1, \tag{2.28}$$

we obtain

$$K_{\epsilon^l}(\omega) = \langle \omega, \delta^{\widehat{\otimes}\epsilon^l} \rangle = \langle \omega, \delta_l \rangle = \langle \omega, \zeta_i(x) p_j(z) \rangle, \tag{2.29}$$

if $l = z(i,j)$.

By Theorem 2.7 any sequence of functions $f_m \in \widehat{L}^2(\pi^{\times m})$, $m = 0, 1, 2, \ldots$, such that $\sum_{m \geq 1} m! \|f_m\|_{L^2(\pi^{\times m})}^2 < \infty$ defines a random variable $F \in L^2(\mu)$ by $F(\omega) = \sum_{m \geq 0} \langle C_m(\omega), f_m \rangle$. Since each $f_m$ is contained in the closure of the linear span of the orthogonal family $\{\delta^{\widehat{\otimes}\alpha}\}_{|\alpha|=m}$ in $L^2(\pi^{\times m})$, we get for all $m \geq 1$ the representation

$$f_m = \sum_{|\alpha|=m} c_\alpha \delta^{\widehat{\otimes}\alpha} \tag{2.30}$$

in $L^2(\pi^{\times m})$ for $c_\alpha \in \mathbb{R}$. Hence, we can restate Theorem 2.7 as follows.

THEOREM 2.8. *The family $\{K_\alpha\}_{\alpha \in \mathcal{J}}$ constitutes an orthogonal basis for $L^2(\mu)$ with norm expression*

$$\|K_\alpha\|_{L^2(\mu)}^2 = \alpha! := \alpha_1! \alpha_2! \cdots, \tag{2.31}$$

*for $\alpha = (\alpha_1, \alpha_2, \ldots) \in \mathcal{J}$. Thus, every $F \in L^2(\mu)$ has the unique representation*

$$F = \sum_{\alpha \in \mathcal{J}} c_\alpha K_\alpha, \tag{2.32}$$

*where $c_\alpha \in \mathbb{R}$ for all $\alpha$ and where we set $c_0 = E[F]$.*

*Moreover, we have the isometry*

$$\|F\|_{L^2(\mu)}^2 = \sum_{\alpha \in \mathcal{J}} \alpha! c_\alpha^2. \tag{2.33}$$

EXAMPLE 2.9. (i) Choose $F(\omega) = \eta(x)$, the $d$-parameter Lévy process. Then $\eta(x) = \int_{[0,x_1] \times \cdots \times [0,x_d] \times \mathbb{R}_0} z \widetilde{N}(dx, dz) = \langle \omega, \chi_{[0,x_1] \times \cdots \times [0,x_d]}(x) \cdot z \rangle$ a.e. and it follows by (2.29) that

$$\eta(x) = \sum_{k \geq 1} m \int_0^{x_d} \cdots \int_0^{x_1} \zeta_k(x_1, \ldots, x_d) \, dx_1 \cdots dx_d \cdot K_{\epsilon^{z(k,1)}}, \tag{2.34}$$

where $m = \|x\|_{L^2(\nu)}$.

(ii) Let $\Lambda_1 \subset \mathbb{R}^m$, $\Lambda_2 \subset \mathbb{R}_0$ with $\pi(\Lambda_1 \times \Lambda_2) < \infty$. Set $f_1(x,z) = \chi_{\Lambda_1 \times \Lambda_2}(x,z)$. Then by (2.29) and (2.30) we get for $F = \widetilde{N}(\Lambda_1, \Lambda_2) = \langle \omega, f_1 \rangle$

$$(2.35) \qquad \widetilde{N}(t, \Lambda) = \sum_{k,m \geq 1} \int_{\Lambda_1} \int_{\Lambda_2} \zeta_k(x) p_m(z) \nu(dz)\, dx \cdot K_{\epsilon^{z(k,m)}}.$$

Next we define various generalized function spaces that relate to $L^2(\mu)$ in a natural way. These spaces turn out to be a useful tool to study stochastic partial differential equations. Our spaces are Lévy versions of the Kondratiev spaces, which were originally introduced in [22]. See also [4] and [25] in the context of Gaussian analysis. The one-parameter case with respect to the Lévy white noise measure $\mu$ can be found in [10] and [34]. The extension to multidimensional parameter sets is analogous.

DEFINITION 2.10. (i) *The stochastic test function spaces.* Let $0 \leq \rho \leq 1$. For an expansion $f = \sum_{\alpha \in \mathcal{J}} c_\alpha K_\alpha \in L^2(\mu)$ define the norm

$$(2.36) \qquad \|f\|_{\rho,k}^2 := \sum_{\alpha \in \mathcal{J}} (\alpha!)^{1+\rho} c_\alpha^2 (2\mathbb{N})^{k\alpha}$$

for $k \in \mathbb{N}_0$, where $(2\mathbb{N})^{k\alpha} = (2 \cdot 1)^{k\alpha_1}(2 \cdot 2)^{k\alpha_2} \cdots (2 \cdot m)^{k\alpha_m}$, if $\mathrm{Index}(\alpha) = m$.

Let

$$(\mathcal{S})_{\rho,k} := \{f : \|f\|_{\rho,k} < \infty\}$$

and define

$$(2.37) \qquad (\mathcal{S})_\rho := \bigcap_{k \in \mathbb{N}_0} (\mathcal{S})_{\rho,k},$$

endowed with the projective topology.

(ii) *The stochastic distribution spaces.* Let $0 \leq \rho \leq 1$. In the same manner, define for a formal expansion $F = \sum_{\alpha \in \mathcal{J}} b_\alpha K_\alpha$ the norms

$$(2.38) \qquad \|F\|_{-\rho,-k}^2 := \sum_{\alpha \in \mathcal{J}} (\alpha!)^{1-\rho} c_\alpha^2 (2\mathbb{N})^{-k\alpha}, \qquad k \in \mathbb{N}_0.$$

Set

$$(\mathcal{S})_{-\rho,-k} := \{F : \|F\|_{-\rho,-k} < \infty\}$$

and define

$$(2.39) \qquad (\mathcal{S})_{-\rho} := \bigcup_{k \in \mathbb{N}_0} (\mathcal{S})_{-\rho,-k},$$

equipped with the inductive topology.

We can regard $(\mathcal{S})_{-\rho}$ as the dual of $(\mathcal{S})_\rho$ by the action

$$\langle F, f \rangle = \sum_{\alpha \in \mathcal{J}} b_\alpha c_\alpha \alpha! \tag{2.40}$$

for $F = \sum_{\alpha \in \mathcal{J}} b_\alpha K_\alpha \in (\mathcal{S})_{-\rho}$ and $f = \sum_{\alpha \in \mathcal{J}} b_\alpha K_\alpha \in (\mathcal{S})_\rho$. Note that for general $0 \leq \rho \leq 1$ we have

$$(\mathcal{S})_1 \subset (\mathcal{S})_\rho \subset (\mathcal{S})_0 \subset L^2(\mu) \subset (\mathcal{S})_{-0} \subset (\mathcal{S})_{-\rho} \subset (\mathcal{S})_{-1}. \tag{2.41}$$

The space $(\mathcal{S}) := (\mathcal{S})_0$, respectively, $(\mathcal{S})^* := (\mathcal{S})_{-0}$, is a Lévy version of the *Hida test function space*, respectively, *Hida stochastic distribution space*. For more information about these or related spaces in the Gaussian and Poissonian case we refer to [16] and [17].

One of the remarkable properties of the space $(\mathcal{S})^*$ is that it accomodates the ($d$-parameter) Lévy white noise. See [10].

DEFINITION 2.11. The ($d$-parameter) *Lévy white noise* $\dot{\eta}(x)$ of the Lévy process $\eta(x)$ (with $m = \|z\|_{L^2(\nu)}$) is defined by the formal expansion

$$\dot{\eta}(x) = m \sum_{k \geq 1} \zeta_k(x) K_{\epsilon^{z(k,1)}}, \tag{2.42}$$

where $\zeta_k(x)$ is defined by Hermite functions, $z(i,j)$ is the map in (2.24) and where $\epsilon^l \in \mathcal{J}$ is defined as in (2.28).

REMARK 2.12. (i) Because of the uniform boundedness of the Hermite functions (see, e.g., [37]) the Lévy white noise $\dot{\eta}(x)$ takes values in $(\mathcal{S})^*$ for all $x$. Further, it follows from (2.34) that

$$\frac{\partial^d}{\partial x_1 \cdots \partial x_d} \eta(x) = \dot{\eta}(x) \qquad \text{in } (\mathcal{S})^*. \tag{2.43}$$

This justifies the name white noise for $\dot{\eta}(x)$.

(ii) Just as in [34] the ($d$-parameter) *white noise* $\dot{\widetilde{N}}(x,z)$ of the Poisson random measure $\widetilde{N}(dx,dz)$ can be defined by

$$\dot{\widetilde{N}}(x,z) = \sum_{k,m \geq 1} \zeta_k(x) p_m(z) \cdot K_{\epsilon^{z(k,m)}}, \tag{2.44}$$

where $p_m(z)$ are the polynomials from (2.23). We have that $\dot{\widetilde{N}}(x,z)$ is contained in $(\mathcal{S})^*$ $\pi$-a.e. The relation (2.35) admits the interpretation of $\dot{\widetilde{N}}(x,z)$ as a Radon–Nikodym derivative, that is, (formally)

$$\dot{\widetilde{N}}(x,z) = \frac{\widetilde{N}(dx,dz)}{dx \times \nu(dz)} \qquad \text{in } (\mathcal{S})^*. \tag{2.45}$$

The last relation entitles us to call $\dot{\tilde{N}}(x,z)$ white noise.

Moreover, $\dot{\eta}(x)$ is related to $\dot{\tilde{N}}(x,z)$ by

$$\dot{\eta}(x) = \int_{\mathbb{R}} z \dot{\tilde{N}}(x,z)\nu(dz). \tag{2.46}$$

The relation above is given in terms of a Bochner integral with respect to $\nu$ (see [34]).

2.3. *Wick product and Hermite transform.* In this section we define a (*stochastic*) *Wick product* on the space $(\mathcal{S})_{-1}$ with respect to the Lévy white noise measure $\mu$. Then we give the definition of the *Hermite transform* and apply it to establish a characterization theorem for the space $(\mathcal{S})_{-1}$.

The Wick product was first introduced by Wick [40] and used as a renormalization technique in quantum field theory. Later on a (stochastic) Wick product was considered by Hida and Ikeda [15]. This subject both in mathematical physics and probability theory is comprehensively treated in Dobroshin and Minlos [11]. Today the Wick product provides a useful concept for a variety of applications, for example, it is important in the study of stochastic ordinary or partial differential equations (see, e.g., [17]).

The next definition is a $d$-parameter version of Definition 3.11 in [10].

DEFINITION 2.13. The *Lévy Wick product* $F \diamond G$ of two elements
$$F = \sum_{\alpha \in \mathcal{J}} a_\alpha K_\alpha, \qquad G = \sum_{\beta \in \mathcal{J}} b_\beta K_\beta \in (\mathcal{S})_{-1} \qquad \text{with } a_\alpha, b_\beta \in \mathbb{R}$$
is defined by
$$F \diamond G = \sum_{\alpha,\beta \in \mathcal{J}} a_\alpha b_\beta K_{\alpha+\beta}. \tag{2.47}$$

REMARK 2.14. Let $f_n = \sum_{|\alpha|=n} c_\alpha \delta^{\widehat{\otimes}\alpha} \in \widehat{L}^2(\pi^{\times n})$ and $g_m = \sum_{|\beta|=m} b_\beta \times \delta^{\widehat{\otimes}\beta} \in \widehat{L}^2(\pi^{\times m})$ according to (2.30). Then we have
$$f_n \widehat{\otimes} g_m = \sum_{|\alpha|=n}\sum_{|\beta|=m} c_\alpha b_\beta \delta^{\widehat{\otimes}(\alpha+\beta)} = \sum_{|\gamma|=n+m}\sum_{\alpha+\beta=\gamma} c_\alpha b_\beta \delta^{\widehat{\otimes}\gamma}$$
in $L^2(\pi^{\times(n+m)})$. Hence,
$$\langle C_n(\omega), f_n \rangle \diamond \langle C_m(\omega), g_m \rangle = \langle C_{n+m}(\omega), f_n \widehat{\otimes} g_m \rangle. \tag{2.48}$$

REMARK 2.15. A remarkable property of the Wick product is that it is implicitly contained in the Itô–Skorohod integrals. The reason for this fact is that if $Y(t) = Y(t,\omega)$ is Skorohod integrable, then (see [10])
$$\int_0^T Y(t)\delta\eta(t) = \int_0^T Y(t) \diamond \dot{\eta}(t)\,dt. \tag{2.49}$$

The left-hand side denotes the Skorohod integral of $Y(t)$ and the integral on the right-hand side is the Bochner integral on $(\mathcal{S})^*$. The Skorohod integral extends the Itô integral in the sense that both integrals coincide, if $Y(t,\omega)$ is adapted, that is, we then have

$$\int_0^T Y(t)\delta\eta(t) = \int_0^T Y(t)\,d\eta(t). \tag{2.50}$$

Note that a version of (2.49) holds for the white noise $\dot{\tilde{N}}(t,x)$, too (see [34]). The extension to the $d$-parameter case is given in [28].

REMARK 2.16. It is important to note that the spaces $(\mathcal{S})_1$, $(\mathcal{S})_{-1}$ and $(\mathcal{S})$, $(\mathcal{S})^*$ form topological algebras with respect to the Lévy Wick product $\diamond$ (for an analogous proof see [35] and [17]). For more information about the Wick product and Skorohod integration in the Poissonian and Gaussian case see, for example, [16, 17] and [31].

The *Hermite transform*, which appeared first in Lindstrøm, Øksendal and Ubøe [27], gives the interpretation of $(\mathcal{S})_{-1}$ in terms of elements in the algebra of power series in infinitely many complex variables. This transform has been applied in many different directions in the Gaussian and Poissonian case (see, e.g., [17]). Its definition for ($d$-parameter) Lévy processes is analogous.

DEFINITION 2.17. Let $F = \sum_{\alpha \in \mathcal{J}} a_\alpha K_\alpha \in (\mathcal{S})_{-1}$ with $a_\alpha \in \mathbb{R}$. Then the *Lévy Hermite transform of* $F$, denoted by $\mathcal{H}F$, is defined by

$$\mathcal{H}F(z) = \sum_{\alpha \in \mathcal{J}} a_\alpha z^\alpha \in \mathbb{C}, \tag{2.51}$$

if convergent, where $z = (z_1, z_2, \ldots) \in \mathbb{C}^{\mathbb{N}}$ (the set of all sequences of complex numbers) and

$$z^\alpha = z_1^{\alpha_1} z_2^{\alpha_2} \cdots z_n^{\alpha_n} \cdots, \tag{2.52}$$

if $\alpha = (\alpha_1, \alpha_2, \ldots) \in \mathcal{J}$, where $z_j^0 = 1$.

EXAMPLE 2.18. We want to determine the Hermite transform of the $d$-parameter Lévy white noise $\dot{\eta}(x)$. Since $\dot{\eta}(x) = m\sum_{k\geq 1} \zeta_k(x) K_{\epsilon^{z(k,1)}}$, we get

$$\mathcal{H}(\dot{\eta})(x,z) = m\sum_{k\geq 1} \zeta_k(x) \cdot z_{z(k,1)}, \tag{2.53}$$

which is convergent for all $z \in (\mathbb{C}^{\mathbb{N}})_c$ (the set of all finite sequences in $\mathbb{C}^{\mathbb{N}}$).

One of the useful properties of the Hermite transform is that it converts the Wick product into ordinary (complex) products.

PROPOSITION 2.19. *If $F$, $G \in (\mathcal{S})_{-1}$, then*

$$\mathcal{H}(F \diamond G)(z) = \mathcal{H}(F)(z) \cdot \mathcal{H}(G)(z) \tag{2.54}$$

*for all $z$ such that $\mathcal{H}(F)(z)$ and $\mathcal{H}(G)(z)$ exist.*

PROOF. The proof is an immediate consequence of Definition 2.13. □

In the following we define for $0 < R, q < \infty$ the infinite-dimensional neighborhoods $K_q(R)$ in $\mathbb{C}^{\mathbb{N}}$ by

$$K_q(R) = \left\{ (\xi_1, \xi_2, \dots) \in \mathbb{C}^{\mathbb{N}} : \sum_{\alpha \neq 0} |\xi^\alpha|^2 (2\mathbb{N})^{q\alpha} < R^2 \right\}. \tag{2.55}$$

By the same proof as in the Gaussian case (see Theorem 2.6.11 in [17]) we deduce the following characterization theorem for the space $(\mathcal{S})_{-1}$.

THEOREM 2.20. (i) *If $F = \sum_{\alpha \in \mathcal{J}} a_\alpha K_\alpha \in (\mathcal{S})_{-1}$, then there are $q, M_q < \infty$ such that*

$$|\mathcal{H}F(z)| \leq \sum_{\alpha \in \mathcal{J}} |a_\alpha| |z^\alpha| \leq M_q \left( \sum_{\alpha \in \mathcal{J}} (2\mathbb{N})^{q\alpha} |z^\alpha|^2 \right)^{1/2} \tag{2.56}$$

*for all $z \in (\mathbb{C}^{\mathbb{N}})_c$.*

*In particular, $\mathcal{H}F$ is a bounded analytic function on $K_q(R)$ for all $R < \infty$.*

(ii) *Conversely, assume that $g(z) = \sum_{\alpha \in \mathcal{J}} b_\alpha z^\alpha$ is a power series of $z \in (\mathbb{C}^{\mathbb{N}})_c$ such that there exist $q < \infty$, $\delta > 0$ with $g(z)$ is absolutely convergent and bounded on $K_q(\delta)$ then there exists a unique $G \in (\mathcal{S})_{-1}$ such that $\mathcal{H}G = g$, namely,*

$$G = \sum_{\alpha \in \mathcal{J}} b_\alpha K_\alpha. \tag{2.57}$$

**3. Application: the stochastic Poisson equation driven by space-time Lévy white noise.** Let us illustrate how the framework, developed in Section 2, can be applied to solve the *stochastic Poisson equation*

$$\begin{aligned} \Delta U(x) &= -\dot{\eta}(x), & x \in D, \\ U(x) &= 0, & x \in \partial D, \end{aligned} \tag{3.1}$$

where $\Delta = \sum_{k=1}^{d} \frac{\partial^2}{\partial x_k^2}$ is the Laplace operator in $\mathbb{R}^d$, $D$ is a bounded domain with regular boundary (see, e.g., Chapter 9 in [33]) and where $\dot{\eta}(x) = m \sum_{k \geq 1} \zeta_k(x) K_{\epsilon^{z(k,1)}}$ is the $d$-parameter Lévy white noise (Definition 2.11).

As mentioned in the Introduction, the model (3.1) gives a description of the temperature $U(x)$ in the region $D$ under the assumption that the temperature at the boundary is kept equal to zero and that there is a white noise heat source in $D$.

Note that $\Delta U(x)$ in (3.1) is defined in the sense of the topology on $(\mathcal{S})_{-1}$.

Now we aim at converting the system (3.1) into a *deterministic* partial differential equation with complex coefficients by applying the Hermite transform (2.51) to both sides of (3.1). Then we try to solve the resulting PDE, and we take the inverse Hermite transform of the solution, if existent, to obtain a solution of the original equation. Before we proceed to realize our strategy, we need the following result.

LEMMA 3.1. *Suppose $X$ and $F$ are functions from $D$ in (3.1) to $(\mathcal{S})_{-1}$ such that*

$$(3.2) \qquad \Delta \mathcal{H} X(x,z) = \mathcal{H} F(x,z)$$

*for all $(x,z) \in D \times K_q(\delta)$ for some $q < \infty$, $\delta > 0$.*

*Furthermore, assume for all $j$ that $\frac{\partial^2}{\partial x_j^2} \mathcal{H} F(x,z)$ is bounded on $D \times K_q(\delta)$, continuous with respect to $x \in D$ for each $z \in K_q(\delta)$ and analytic with respect to $z \in K_q(\delta)$ for all $x \in D$.*

*Then*

$$(3.3) \qquad \Delta X(x) = F(x) \qquad \text{for all } x \in D.$$

PROOF. Use repeatedly the same proof of Lemma 2.8.4 in [17] in the case of higher-order derivatives. □

Now, we take the Hermite transform of (3.1) and we get

$$(3.4) \qquad \begin{aligned} \Delta u(x,z) &= -\mathcal{H}(\dot{\eta})(x,z), & x \in D, \\ u(x,z) &= 0, & x \in \partial D, \end{aligned}$$

where $u = \mathcal{H} U$ and $\mathcal{H}(\dot{\eta})(x,z) = m \sum_{k \geq 1} \zeta_k(x) \cdot z_{z(k,1)}$ for $z \in (\mathbb{C}^{\mathbb{N}})_c$ [see (2.53)]. By comparing the real and imaginary parts of equation (3.4), one checks that

$$(3.5) \qquad u(x,z) = \int_{\mathbb{R}^d} G(x,y) \cdot \mathcal{H}(\dot{\eta})(y,z) \, dy,$$

where $G(x,y)$ is the classical Green function of $D$ with $G = 0$ outside of $D$ (see, e.g., Chapter 9 in [33]). Since $G(x,\cdot) \in L^1(\mathbb{R}^d)$ for all $x$, the right-hand side of (3.5) exists for all $z \in (\mathbb{C}^{\mathbb{N}})_c$ and $x \in D$. Hence, $u(x,z)$ is defined for such $z$, $x$.

Further, we see for all $z \in (\mathbb{C}^{\mathbb{N}})_c$ that

$$
\begin{aligned}
|u(x,z)| &\leq \sum_k |z_{z(k,1)}| \int_{\mathbb{R}^d} |G(x,y)||\zeta_k(y)|\,dy \\
&\leq \text{const.} \sum_k |z^{\epsilon_k}| \\
&\leq \text{const.} \left(\sum_k |z^{\epsilon_k}|^2 (2\mathbb{N})^{2\epsilon_k}\right)^{1/2} \left(\sum_k (2\mathbb{N})^{-2\epsilon_k}\right)^{1/2} \\
&\leq \text{const.} \cdot R \cdot \left(\sum_k (2k)^{-2}\right)^{1/2} < \infty
\end{aligned}
\tag{3.6}
$$

for all $z \in K_2(R)$. Besides this (3.5) shows that $u(x,z)$ is analytical in $z$. Thus, we conclude by the characterization theorem (Theorem 3.8) that there exists a function $U : \overline{D} \to (\mathcal{S})_{-1}$ such that $\mathcal{H}U(x,z) = u(x,z)$. Next we want to verify the assumptions of Lemma 3.1 for $X = U$ and $F = -\dot{\eta}$. It is known from the general theory of deterministic elliptic PDE's (see, e.g., [7]) that for all open and relatively compact $V$ in $D$ there exists a $C$ such that

$$
\|u(\cdot,z)\|_{C^{2+\alpha}(V)} \leq C(\|\Delta u(\cdot,z)\|_{C^{\alpha}(V)} + \|u(\cdot,z)\|_{C(V)})
\tag{3.7}
$$

for all $z \in (\mathbb{C}^{\mathbb{N}})_c$. Since $\Delta u = -\mathcal{H}\dot{\eta}$ and $u$ are bounded on $D \times K_2(R)$, it follows that $\frac{\partial^2}{\partial x_j^2} u(x,z)$ is bounded for such $x$, $z$. Thus, by Lemma 3.1 $U$ is a solution of system (3.1).

Further, we follow from Lemma 3.18 in [10] that the Bochner integral $\int_{\mathbb{R}^d} G(x,y)\dot{\eta}(x)\,dx$ exists in $(\mathcal{S})^*$ (see Definition 3.16 in [10]) and that

$$
\int_{\mathbb{R}^d} G(x,y)\dot{\eta}(y)\,dy = m \sum_{k \geq 1} \int_{\mathbb{R}^d} G(x,y)\zeta_k(y)\,dy\, K_{z(k,1)}.
\tag{3.8}
$$

Then one realizes that the right-hand side of (3.5) is the Hermite transform of (3.8).

So we obtain the following result.

THEOREM 3.2. *There exists a unique stochastic distribution process* $U : \overline{D} \to (\mathcal{S})^*$, *solving system* (3.1). *The solution is twice continuously differentiable in* $(\mathcal{S})^*$ *and takes the form*

$$
U(x) = \int_{\mathbb{R}^d} G(x,y)\dot{\eta}(y)\,dy = m \sum_{k \geq 1} \int_{\mathbb{R}^d} G(x,y)\zeta_k(y)\,dy\, K_{z(k,1)},
\tag{3.9}
$$

*where* $m = \|z\|_{L^2(\nu)}$.

We conclude this section with a remark about an alternative approach to SPDEs driven by Lévy space-time white noise.

REMARK 3.3. Let us briefly describe how the concepts in [28] can be used to establish a framework similar to Section 2. Instead of the spaces $(\mathcal{S})_{-\rho}$, consider the distribution spaces in [28] and instead of the $\mathcal{H}$-transform, use the $\mathcal{S}$-transform in [28]. The $\mathcal{S}$-transform, is of the form

$$\mathcal{S}(F)(\phi) = \langle\!\langle F(\omega), \widetilde{e}(\phi, \omega)\rangle\!\rangle$$

for distributions $F$ and for $\phi$ in an open neighborhood of zero in $\widetilde{\mathcal{S}}(X)$, where the function $\widetilde{e}(\phi, \omega)$ is as in (2.11) and where $\langle\!\langle \cdot, \cdot \rangle\!\rangle$ is an extension of the innerproduct on $L^2(\mu)$. Moreover, the process $\dot{\eta}(x)$ can be replaced by

$$\dot{\eta}(x) := \langle C_1(\omega), z\delta_x \rangle$$

and the white noise $\dot{\widetilde{N}}$ can be defined by

$$\dot{\widetilde{N}}(x, z) := \langle C_1(\omega), \delta_{(x,z)} \rangle,$$

where $\delta_y$ is the Dirac measure in a point $y$. Further, by the properties of the $\mathcal{S}$-transform (see [28]) one can prove a similar result as Lemma 3.1. Moreover, the $\mathcal{S}$-transform of $\dot{\eta}(x) := \langle C_1(\omega), z\delta_x \rangle$ is

$$\mathcal{S}(\dot{\eta}(x))(\phi) = \int_{\mathbb{R}_0} \phi(x, z)\nu(dz)$$

(see proof of Proposition 7.5 in [28]). Hence, we can solve system (3.1) by finding a function $u$ such that

$$\begin{aligned}\Delta u(x, \phi) &= -\int_{\mathbb{R}_0} \phi(x, z)\nu(dz), & x \in D, \\ u(x, \phi) &= 0, & x \in \partial D.\end{aligned}$$

The obvious candidate for $u$ is given by the Green function $G$:

$$u(x, \phi) = \int_{\mathbb{R}^d} G(x, y) \int_{\mathbb{R}_0} \phi(x, z)\nu(dz)\, dy.$$

Hence, the solution is given by the inverse $\mathcal{S}$-transform, yielding the same result as in Theorem 3.2 for all Lévy measures. Moreover, within a similar setting one can solve more general versions of the problem. However, the use of the $\mathcal{H}$-transform has some advantages. For instance, it enables the application of methods of complex analysis.

**4. Discussion of the solution.** As mentioned in the Introduction, our interpretation of the solution $U(x) = U(x,\omega)$ of (3.1) is the following:

For each $x$ we have $U(x,\cdot) \in (\mathcal{S})^*$ and $x \mapsto U(x)$ satisfies (3.1) in the strong sense as an $(\mathcal{S})^*$-valued function.

Regarding the solution itself, given by (3.9), the interpretation is the following:

For each $x$, $U(x)$ is a stochastic distribution whose action on a stochastic test function $f \in (\mathcal{S})$ is

$$\langle U(x), f \rangle = \int_{\mathbb{R}^d} G(x,y) \langle \dot{\eta}(y), f \rangle \, dy, \tag{4.1}$$

where

$$\langle \dot{\eta}(y), f \rangle = \left\langle m \sum_{k \geq 1} \xi_k(y) K_{\epsilon^{z(k,1)}}, f \right\rangle$$

$$= m \sum_{k \geq 1} \xi_k(y) E[K_{\epsilon^{z(k,1)}} f] \qquad [\text{see (2.40)}].$$

Relation (4.1) gives rise to the intepretation that the solution $U(x,\omega)$ takes $\omega$-averages for all $x$.

In general, we are not able to represent this stochastic distribution as a classical random variable $U(x,\omega)$. However, if the space dimension $d$ is low, we can say more:

COROLLARY 4.1. *Suppose $d \leq 3$. Then the solution $U(x,\cdot)$ given by (3.9) in Theorem 3.2 belongs to $L^2(\mu)$ for all $x$ and is continuous in $x$.*

*Moreover,*

$$U(x) = \int_{\mathbb{R}^d} G(x,y) \, d\eta(y). \tag{4.2}$$

PROOF. Since the singularity of $G(x,y)$ at $y = x$ has the order of magnitude $|x-y|^{2-d}$ for $d \geq 3$ and $\log \frac{1}{|x-y|}$ for $d = 2$ (with no singularity for $d = 1$), we see by using polar coordinates that

$$\int_D G^2(x,y) \, dy \leq C \int_0^1 r^{2(2-d)} r^{d-1} \, dr$$

$$= \int_0^1 r^{3-d} \, dr$$

$$< \infty$$

for $d \leq 3$. Hence, by Remark 2.15 we get

$$U(x) = \int_{\mathbb{R}^d} G(x,y) \diamond \dot{\eta}(y) \, dy = \int_{\mathbb{R}^d} G(x,y) \, d\eta(y)$$

and

$$E[U^2(x)] = E\left[\left(\int_{\mathbb{R}^d} G(x,y)\,d\eta(y)\right)^2\right] = M \int_D G^2(x,y)\,dy < \infty$$

for $d \leq 3$. $\square$

Some remaining natural questions are the following:

Q1. If $d \geq 4$, is $U(x) \in L^p(\mu)$ for some $p = p(d) \geq 1$?

Q2. Is it possible to prove that equation (3.3) also has a solution $\widetilde{U}(x) = \widetilde{U}(x,\omega)$ of Walsh type, that is, such that, for some $n$,

$$x \mapsto \widetilde{U}(x,\omega) \in H^{-n}(D) \qquad \text{for a.a. } \omega$$

and $\widetilde{U}(x)$ solves (3.1) in the (classical) sense of distributions, for a.a. $\omega$? [See (1.3) in the analogous Brownian motion case.]

We will not pursue any of these questions here.

**Acknowledgments.** We wish to thank the anonymous referees whose constructive criticism has been particularly appreciated. We thank F. E. Benth for suggestions and valuable comments. We are also grateful to G. Di Nunno for helpful remarks.

A. LØKKA
F. PROSKE
DEPARTMENT OF MATHEMATICS
UNIVERSITY OF OSLO
P.O. BOX 1053 BLINDERN
N-0316 OSLO
NORWAY
E-MAIL: alokka@math.uio.no
E-MAIL: proske@math.uio.no

B. ØKSENDAL
DEPARTMENT OF MATHEMATICS
UNIVERSITY OF OSLO
P.O. BOX 1053 BLINDERN
N-0316 OSLO
NORWAY
AND
NORWEGIAN SCHOOL OF ECONOMICS
AND BUSINESS ADMINISTRATION
HELLEVEIEN 30
N-5045 BERGEN
NORWAY
E-MAIL: oksendal@math.uio.no